\documentclass[12pt]{amsart}
\newtheorem{theorem}{Theorem}[section]

\newcommand{\ncom}{\newcommand}
\ncom{\lrar}{\longrightarrow}
\ncom{\rar}{\rightarrow}
\ncom{\ov}{\overline}
\ncom{\m}{\mbox}
\ncom{\sta}{\stackrel}
\ncom{\comx}{{\mathbb C}}
\ncom{\A}{{\mathbb A}}
\ncom{\Z}{{\mathbb Z}}
\ncom{\Q}{{\mathbb Q}}
\ncom{\R}{{\mathbb R}}
\ncom{\G}{{\mathbb G}}
\ncom{\al}{\alpha}
\ncom{\p}{{\mathbb P}}
\ncom{\E}{{\mathbb E}}
\ncom{\N}{{\mathbb N}}
\ncom{\K}{{\mathbb K}}
\ncom{\X}{{\mathbb X}}
\ncom{\f}{\frac}
\ncom{\cA}{{\mathcal A}}
\ncom{\cB}{{\mathcal B}}
\ncom{\cD}{{\mathcal D}}
\ncom{\cX}{{\mathcal X}}
\ncom{\cO}{{\mathcal O}}
\ncom{\cW}{{\mathcal W}}
\ncom{\cL}{{\mathcal L}}
\ncom{\cP}{{\mathcal P}}
\ncom{\cH}{{\mathcal H}}
\ncom{\cS}{{\mathcal S}}
\ncom{\cM}{{\mathcal M}}
\ncom{\cC}{{\mathcal C}}
\ncom{\cT}{{\mathcal T}}
\ncom{\cF}{{\mathcal F}}
\ncom{\cN}{{\mathcal N}}
\ncom{\cJ}{{\mathcal J}}
\ncom{\cK}{{\mathcal K}}
\ncom{\cV}{{\mathcal V}}
\ncom{\cZ}{{\mathcal Z}}
\ncom{\cU}{{\mathcal U}}
\ncom{\cSU}{{\mathcal S \mathcal U}}
\ncom{\cG}{{\mathcal G}}
\ncom{\cQ}{{\mathcal Q}}
\ncom{\cR}{{\mathcal R}}
\ncom{\cY}{{\mathcal Y}}
\ncom{\cE}{{\mathcal E}}
\ncom{\what}{\widehat}
\ncom{\delbar}{\overline{\partial}}

\ncom{\eop}{{\hfill $\Box$}}

\begin{document}
\baselineskip=16pt

\title{Regulators of canonical extensions are torsion: the smooth divisor case}


\author[J. N. Iyer]{Jaya NN  Iyer}

\address{The Institute of Mathematical Sciences, CIT
Campus, Taramani, Chennai 600113, India}
\email{jniyer@imsc.res.in}

\footnotetext{Mathematics Classification Number: 14C25, 14D05, 14D20, 14D21 }

\begin{abstract}
In this note, we report on a work jointly done with C. Simpson on a generalization of Reznikov's theorem which says that 
the Chern-Simons classes and in particular the Deligne Chern classes (in degrees $>\,1$) are torsion, of a 
flat vector bundle on a smooth complex projective variety. We consider the case of a smooth quasi--projective variety 
with an irreducible smooth divisor at infinity. We define
the Chern-Simons classes of the Deligne's \textit{canonical extension} of a flat vector bundle with unipotent monodromy 
at infinity, which lift the Deligne Chern classes and prove that these classes are torsion. The details of the proof can be found in arxiv: math.AG.07070372. 
\end{abstract}
\maketitle

\section{Introduction}


This is a report on a work jointly done with C. Simpson; ``Regulators of canonical extensions are torsion: the smooth divisor case'', arXiv:math.AG.07070372, at IAS, Princeton during the special year in Algebraic geometry 2006-07. This is a continuation of \cite{IySi}, \cite{IySi2}, motivated by Reznikov's theorem \cite{Reznikov}  which shows the triviality of the Chern-Simons classes of flat vector bundles on complex smooth projective varieties in the $\R/\Q$-cohomology. This answers a question of Cheeger-Simons in the projective case and in turn proves a conjecture of Bloch on the triviality of the Chern classes in the rational Deligne cohomology. We would like to conclude the same for Deligne's canonical extension with unipotent local monodromies around a normal crossing divisor at infinity. We treat the case when we have a smooth divisor at infinity. 

These questions can be historically traced back to Weil's theorem \cite[p.57]{Chern}.
This theorem states that the de Rham Chern classes $c_i^{dR}$ which are obtained by substitution of the curvature form of a connection $\theta$ in the $GL_n$-invariant polynomials $P_i$ (a construction due to Chern \cite{Chern}), of a complex vector bundle $E$ are independent of the connection. In particular, this says that if $E$ has a flat connection then the de Rham Chern classes are zero and via the de Rham isomorphism the Chern classes in the integral Betti cohomology are torsion. 
 
Since then attempts to  construct secondary invariants for bundles with connection has been of wide interest. A first construction in this direction was done by Chern and Simons \cite{ChernSimons}. These are differential forms, denoted by $TP_i(\theta)$, which live on the total space of $E$. Later in \cite{CheegerSimons},
Cheeger and Simons defined
differential characters $\what{c_i}(E,\theta)$ which lived on the base manifold and which are closely related to the forms $TP_i(\theta)$.  
If $\theta$ is flat then these differential characters actually lie in 
the odd degree $\R/\Z$-cohomology of the base manifold and are called as the Chern-Simons
classes. These are uniquely determined lifts of the Betti Chern class related via the coefficient sequence $0\rar \Z\rar \R \rar \R/\Z\rar 0$. A statement similar to Weil's theorem holds; the Chern-Simon's classes
in degrees at least two are constant in a family of flat connections \cite{CheegerSimons}. This is usually called as the rigidity property. Since there are countably many connected 
components of the space of representations of the fundamental group of the manifold, the following question was posed by Cheeger-Simons:

\textit{Question}: Suppose $(E,\theta)$ is a flat vector bundle on a smooth manifold $M$.
Are the Chern-Simons classes $\what{c_i}(E,\theta)\,\in\, H^{2i-1}(M,\R/\Z)$
torsion, when $i\geq 2$.

Our aim here is firstly to extend this question when $X$ is smooth and quasi--projective with an irreducible smooth divisor $D$ at infinity. We consider a flat bundle on $X-D$ which has unipotent monodromy around the divisor $D$. We define the Chern-Simons classes of the canonical extension \cite{De} on $X$ which extend the classes on $X-D$ of the flat connection. Furthermore, the extended classes are shown to lift the Deligne Chern classes whenever $X$ is projective.

Our main theorem is
\begin{theorem}\label{maintheorem}
Suppose $X$ is a smooth quasi--projective variety defined over $\comx$. Let $(E,\nabla)$ be a flat connection on $U:=X-D$ associated to a representation $\rho:\pi_1(U)\rar GL_r(\comx)$. Assume that $D$ is a smooth and irreducible divisor and $(\ov{E},\ov{\nabla})$ be the Deligne canonical extension on $X$ with unipotent monodromy around $D$. Then 
the Chern-Simons classes
$$
\what{c_p}(\rho / X )\,\in\,H^{2p-1}(X , \R / \Z )
$$    
of $(\ov{E},\ov{\nabla})$ are torsion, for $p>1$. If, furthermore, $X$ is projective then 
the Chern classes of $\ov E$ are torsion in
the Deligne cohomology of $X$, in degrees $>1$.
\end{theorem}

What we do here can easily be generalized to the case when $D$ is smooth and has several
disjoint irreducible components.

After posting the preprint on the math arXiv, Deligne \cite{Deligneletter} has indicated a construction of the Chern-Simons classes in general and Esnault has  independently given a construction in \cite{Esnault}.


\section{Proof of Theorem \ref{maintheorem}}

The idea of the proof is to adapt Reznikov's proof for the canonical extension.
For this purpose, we briefly recall his proof.

\subsection{Reznikov's theorem}\cite{Reznikov}
Suppose $X$ is a complex smooth projective variety and $(E,\nabla)$ is a flat connection on $X$. Let $\rho$ be the monodromy representation of the flat connection. 
The main steps in his proof are as follows:

$\bullet$ there is a classifying map
$$
X \sta{\psi_\rho}{\lrar}BGL(\comx)
$$
such that $\what{c_i}(E,\nabla)=\psi_\rho^*\what{c_i}^{univ}$.
Here $GL(\comx)$ is endowed with the discrete topology and $\what{c_i}^{univ}$ are the Beilinson's universal classes.

$\bullet$ it suffices to look at $SL_n$-valued representations and via the rigidity property that it is defined over a number field $F$. So we can assume that
$\rho$ gives a map
$$
X \sta{\psi_\rho}{\lrar}BSL(F).
$$
  
$\bullet$ apply Borel's theorem to say that the real cohomology $H^*(BSL(F),\R)$ is generated by $\sigma^*Vol_{2i-1}$, for $i\geq 1$, and for various embeddings
$F\sta{\sigma}{\hookrightarrow} \comx$.
Here $Vol_{2i-1}$ are the Borel's volume regulators.

$\bullet$ apply Simpson's theorem \cite{Si2} to deform $\rho$ to a complex VHS
and again by the rigidity property, it suffices to assume that $\psi_\rho$ factors via the map
$$
X {\lrar}BU(F).
$$
Here we use the fact that the monodromy representation of a complex VHS takes values in the unitary
group $U_{p,q}$. Observe that $H^{2i-1}(BU(\comx),\R)=0$, for $i\geq 1$.

All the above facts put together imply that the pullback homomorphism
$$
H^*(BSL(F),\Q)\rar H^*(X,\Q)
$$
is the zero map. This suffices to conclude that the Chern-Simons classes of $(E,\nabla)$ are torsion, in degrees at least two.

\subsection{Sketch of proof of Theorem \ref{maintheorem}}

Since the Chern-Simons classes are topological classes, we want to view $(X,D)$
topologically and see the topological information that the datum of the canonical extension gives.

We will consider the following situation which will help us to do similar constructions with the classifying spaces in \textbf{Step 4} below. Suppose $X$ is a smooth manifold, and 
$D\subset X$ is a connected smooth closed subset of real codimension $2$. Let $U:= X-D$ and suppose we can choose a reasonable
tubular neighborhood $B$ of $D$. 
Let $B^{\ast} := B\cap U = B-D$.
It follows that $\pi _1(B^{\ast})\rightarrow \pi _1(B)$ is surjective.
 The diagram
\begin{equation}
\label{pushout}
\begin{array}{ccc}
B^{\ast} & \rightarrow & B \\
\downarrow & & \downarrow \\
U & \rightarrow & X
\end{array}
\end{equation}
is a homotopy pushout diagram. Note also that $B$ retracts to $D$, and $B^{\ast}$ has a tubular structure: 
$$
B^{\ast} \cong S \times (0,1)
$$
where $S\cong \partial B$ is a circle bundle over $D$.

We say that $(X,D)$ is {\em complex algebraic} if $X$ is a smooth complex quasiprojective variety and $D$ an irreducible smooth divisor. 

\textbf{Step 1} (canonical extension on $X$):

Suppose we are given a representation $\rho : \pi _1(U)\rightarrow GL_r(\comx )$, corresponding to a local system $L$ over $U$,
or equivalently to a vector bundle with flat connection $(E,\nabla )$. Let $\gamma$ be a loop going out from the basepoint to a point near $D$,
once around, and back. Then $\pi _1(B)$ is obtained from $\pi _1(B^{\ast})$ by adding the relation $\gamma \sim 1$. 
We assume that {\em the monodromy of $\rho$ at infinity is unipotent}, by which we mean that $\rho (\gamma )$ should be unipotent.

In this situation, there is a canonical and natural way to extend the bundle $E$ to a bundle $\ov E$ over $X$, known as the 
{\em Deligne canonical extension} \cite{De}. The connection $\nabla$ extends
to a connection $\ov{\nabla}$ whose singular terms involved look locally like $Nd\theta$ where $\theta$ is the angular coordinate around $D$. 
In an appropriate frame the singularities of $\ov \nabla $ are only in the strict upper triangular region of the
connection matrix.  In the complex algebraic case, $(E,\nabla )$ are holomorphic, and indeed algebraic with algebraic structure
uniquely determined by the requirement that $\nabla$ have regular singularities. The extended bundle $\ov E$ is
algebraic on $X$ and
$\ov\nabla$ becomes a logarithmic connection \cite{De}. 

\textbf{Step 2} (defining extended regulator classes via patched connection):

We will define {\em extended regulator classes} 
$$
\what{c}_p(\rho / X) \in H^{2p-1}(X, \comx / \Z )
$$
which restrict to the usual regulator classes on $U$. Their imaginary parts define {\em extended volume regulators}
which we write as $Vol_{2p-1}(\rho / X)\in H^{2p-1}(X, \R )$. 

The technique for defining the extended regulator classes is to construct a {\em patched connection} $\nabla ^{\#}$ over $X$.
This will be a smooth connection, however it is not flat. Still, the curvature comes from the singularities of 
$\ov \nabla$ which have been smoothed out, so the curvature is upper-triangular. In particular, the Chern forms for 
$\nabla ^{\#}$ are still identically zero. The Cheeger-Simons theory of differential characters provides a
class of $\nabla ^{\#}$ in the group of differential characters, mapping to the group of closed forms. Since the image,
which is the Chern form, vanishes, the differential character lies in the kernel of this map which is exactly
$H^{2p-1}(X, \comx / \Z )$ \cite[Cor. 2.4]{CheegerSimons}. This is the construction of the regulator class.

\textbf{Step 3} (the extended regulator class lift the Deligne Chern class)

The proof of Dupont-Hain-Zucker that the regulator class lifts the Deligne Chern class, goes through word for word here
to show that this extended regulator class lifts the Deligne Chern class of the canonical extension $\overline{E}$
in the complex algebraic case. For this part, we need $X$ projective. 

\textbf{Step 4} (extended regulator class via $K$-theory):

We also give a different construction of the regulator classes, using the deformation theorem in $K$-theory.
The filtration which we will use to define the patched connection, also leads to a 
polynomial deformation on $B^{\ast}$ between
the representation $\rho$ and its associated-graded. Then, using the fact that $BGL(F[t])^+$ is homotopy-equivalent
to $BGL(F)^+$ and the fact that the square \eqref{pushout} is a homotopy pushout, this 
allows us to construct a map from 
$X$ 
to the homotopy pushout space $BGL(F)^+_{def}$,
$$
X\sta{\psi_\rho}{\rar}BGL(F)^+_{def}.
$$
 A deformation theorem in $K$-theory allows us to identify the cohomology of the pushout space with that of $BGL(F)^+$ and hence we can pull back the universal regulator classes via the map $\psi_\rho$. 
Corollary $7.5$ in \cite{IySi3} says that these are the same as the extended
regulators defined by the patched connection.  

\textbf{Step 5} (rigidity property and deformation to a complex VHS):
 We apply Mochizuki's theorem that any representation can be deformed to a complex variation of Hodge structure, in the quasi-projective case
\cite{Mochizuki}.
The counterpart of the deformation
construction in hermitian $K$-theory
allows us to conclude that the extended volume regulator is zero whenever $\rho$ underlies a complex variation of Hodge structure
in the complex algebraic case.
 This uses the one-variable nilpotent and $SL_2$-orbit theorems and a polynomial deformation as in \textbf{Step 4}, for constructing a map from $X$ to the pushout space $BO(\comx)^+_{def}$ factoring the map $\psi_{\rho}$, and using Karoubi's deformation theorem to identify the cohomologies of $BO(\comx)^+$ with that of  $BO(\comx)^+_{def}$. 
A rigidity statement for the patched connections is discussed and proved in more generality in \cite[\S 6]{IySi3}.

All of the ingredients of Reznikov's original proof \cite{Reznikov} are now present for the extended classes and we can show that the extended regulator classes $\what{c}_p(\rho / X)$ are torsion, in degrees at least two.

Thus we show the generalization of Reznikov's result, in the single divisor case.


\end{document}